\documentclass[11pt]{article}
\usepackage{authblk} % for author information

\usepackage[latin1]{inputenc}
\usepackage[T1]{fontenc}
\usepackage{amsmath}
\usepackage{amsfonts}
\usepackage{amssymb}

\usepackage[normalem]{ulem} % delete underlines in references

\usepackage{amsmath,xcolor,ulem}
\usepackage{amsfonts}
\usepackage{amssymb,todonotes}
\usepackage{amsthm}
\usepackage{graphicx}
\usepackage{marginnote}
\usepackage{float}
\usepackage{subcaption}
\usepackage{xcolor}

\usepackage[colorlinks=true, allcolors=blue]{hyperref}

\usepackage[top=2.8cm,bottom=2.8cm,left=2.5cm,right=2.5cm]{geometry}
\usepackage{float}

\usepackage[font=footnotesize,width=\linewidth]{caption} 
\usepackage{algorithm}
\usepackage{algpseudocode}

\newcommand{\R}{\mathbb R}

\DeclareMathOperator{\GP}{\text{GP}}

\newtheorem{theorem}{Theorem}
\newtheorem{definition}[theorem]{Definition}

\newtheorem*{remark}{Remark}

\title{A Consensus-based optimization algorithm using Gaussian processes for global optimization problems in Sobolev spaces}

\author{Mahmoud Khatab\footnote{University of Wuppertal, Germany, \href{mailto:khatab@uni-wuppertal.de}{khatab@uni-wuppertal.de}}  \; and \;
Claudia Totzeck\footnote{University of Wuppertal, Germany, \href{mailto:totzeck@uni-wuppertal.de}{totzeck@uni-wuppertal.de}}}
\begin{document}
\maketitle

\begin{abstract} 
% To be written at the end...
We propose an algorithm to approximate solutions of global optimization problems in Sobolev spaces that follows the spirit of Consensus-based algorithms in finite dimensions. The main ingredient are Gaussian processes. In fact, we exploit their rich toolbox in order to draw sample functions from Sobolev spaces that satisfy initial values, boundary conditions or state constraints. Well-known marginalization properties of Gaussian processes help us to discretize the algorithm, that is stated in infinite dimensions, appropriately. We illustrate the performance of the algorithm and show its feasibility for nonlinear boundary value problems with state constraints as well as nonlinear optimal control problems constrained by a system of ordinary differential equations with several numerical results.
\end{abstract}

\section{Introduction}
About one decade ago the Consensus-based optimization (CBO) method for nonconvex continuous optimization problems in $\R^d$ was proposed in \cite{cbo2,cbo}. Since then the method has attracted many researchers from different backgrounds to push the limits in the theory of provably convergent methods for global optimization and to extend the range of application of such methods. Let us mention some recent contributions that originate from various mathematical disciplines, i.e., the analysis of partial differential equations \cite{bonandin2025strong,fornasier2025regularity,gerber2025uniform,huang2025faithful}, in particular mean-field theory \cite{baldi2025alternative,gerber2025mean,koss2025mean}, stochastic analysis \cite{bellavia2025discrete,goettlich2025exponential} or mathematical machine learning \cite{carrillo2024fedcbo}. Also the range of applications that may be tackled with CBO-based techniques was extended, see \cite{bungert2025mirrorcbo,bungert2025polarized,chenchene2025consensus,de2025consensus,fornasier2025pde,sun2026consensus}, and more robust versions were proposed \cite{garcia2025defending,huang2026fast}. For a quick start into coding the method we recommend the software package~\cite{Bailo2024joss}. 
Although the approach presented below uses different techniques, the article aligns with recent advances towards global optimization in infinite dimensions. In particular, we want to mention:  the multi-grid finite element approach combined with CBO for the constrained $p$-Allen-Cahn problem with obstacle \cite{beddrich2026}; the CBO approach for global optimization in the space of Gaussian probability measures \cite{borghi2026variational}; the multiscale approach for multi-level optimization with CBO \cite{herty2025multiscale}; and the model predictive control approach enhanced by CBO  \cite{Borghi2025cboMPC}.
%%%%

The main contribution of the present article is the extension of Consensus-based optimization in the direction of global optimization problems in function spaces, that is, to infinite dimensional problems. In particular, we consider Sobolev Spaces with elements that may be represented by samples from Gaussian processes. We remark that in general, the well-posedness analysis of optimization problems in function spaces exploits properties such as coercivity and weakly lower semicontinuity of the objective functional and weak compactness arguments \cite{Hinze2009}. The proofs are therefore not constructive and first-order optimality conditions are typically used to characterize solutions. Moreover, the first-order optimality conditions are the foundation to derive first-order optimization methods such as gradient descent. These are well-understood, especially in the case of Hilbert spaces, and known to globally converge to stationary points. However, additional assumptions, e.g.~convexity, are required to deduce that stationary points are global minimizers.

With the approach proposed here we follow a different avenue, instead of exploiting first-order conditions, we increase the number of function evaluations and propose a zero-order method that aims to approximate the global minimum directly. In fact, we use Gaussian processes to randomly generate Sobolev functions and therefore to translate the exploration and exploitation of Consensus-based dynamics from finite dimensions to the Sobolev case. This allows us to formulate a algorithm directly in the function space. Moreover, the well-established toolbox for Gaussian processes can be exploited to derive easy-to-implement algorithms as we will show with the help of numerical simulation results for ODE- and PDE-boundary value problems as well as optimal control problems constrained by nonlinear interacting particle systems. Before jumping into the details, we state the problem considered  and sketch the main idea that is further explained in the following sections.

\subsection{Problem formulation} 
Let $\mathcal U$ be a function space with elements that can be represented as sample path of measurable Gaussian process (GP) and $f\colon \mathcal U \rightarrow \mathbb R$ a possibly nonconvex functional with unique global minimizer in $\mathcal U$.  We are interested in the global optimization problem
\begin{equation}\label{eq:globalProblem}
\min\limits_{u\in \mathcal U} f(u).
\end{equation}
For instance, $\mathcal U$ can be a Sobolev space of integer order $W^{m,p}(\mathcal D)$, $m\in\mathbb N_0$, $1<p<+\infty$, where $\mathcal D$ is an arbitrary open set \cite{Henderson2024}. In particular, we are interested in the case $p=2$, where the Hilbert structure allows for a direct relation of functions in the reproducing kernel Hilbert space $H_k$ of a given kernel $k$ and their equivalence class in $H^m(\mathcal D)$. The key element, which relates the state space $\mathcal U$ of the optimization problem with the GP is the kernel $k$ also known as covariance function.

\subsection{Sketch of the main idea}
In contrast to the usual use of GP in machine learning \cite{Rasmussen2006}, GP are used in the following to sample functions with required smoothness to lift both, the diffusive part of the CBO dynamics and the initial conditions to the function space setting. Indeed, let $\GP_1$ and $\GP_2$ be Gaussian processes with specific properties that we define later. Further, let the time step size $\tau>0$ and drift parameter $\lambda>0$ be given. To approximate the global solution of \eqref{eq:globalProblem}, we propose the time-discrete dynamics
\begin{subequations} \label{eq:timediscrete}
\begin{align} 
&dU_{j+1}^i = -\lambda (U_j^i - v_f^\alpha(U_j))\tau + \sqrt{2\tau} \|U_j^i - v_f^\alpha(U_j)\|_{\mathcal U} \,\xi_j^i, \\
&U_0^i \sim \GP_1, \quad \xi_j^i \sim \GP_2
, \quad i=1,\dots,N, \quad j=0,\dots,J,
\end{align}
\end{subequations}
where $v_f^\alpha(U_j)$ is the weighted mean, computed in the usual spirit of Consensus-based methods \cite{cbo}
$$ v_f^\alpha(U_j) = \frac{\sum_{i=1}^N U_j^i e^{-\alpha f(U_j^i)}}{\sum_{i=1}^N e^{-\alpha f(U_j^i)}}, \qquad \alpha \gg 0.$$
The similarity to the dynamics in the finite dimensional setting is obvious. The key idea leading to the formulation of the dynamics in \eqref{eq:timediscrete} is strongly related to the interpretation of a GP as a distribution over functions. Indeed, we expect $\xi_j^i$ to mimic the Brownian increments of the finite dimensional scheme. The two Gaussian processes $\GP_1$ and $\GP_2$ play a crucial role and have to be tailored for each state space $\mathcal U$, that is each application individually. To gain further insights, we recall some characteristics of the GP in the next section.

\section{Compilation of properties of Gaussian Processes}
In this section we recall some definitions in the context of GP and related properties that are crucial for the proposed algorithm. From the theoretical point of view a particular focus lies on the relationship of the kernel function and the smoothness of the sample path, and from a numerical point of view on the marginalization property. The latter plays a crucial role for the space discretization of \eqref{eq:timediscrete} and is heavily exploited in the implementation of our algorithm. We closely follow \cite{Henderson2024}.

\subsection{Basic definitions}
Let $(\Omega,\mathcal F, \mathbb P)$ denote a probability space. 
 \begin{definition}\label{def:law}
 The law $\mathbb P_X$ of a random variable $X \colon \Omega \rightarrow \mathbb R$ is the pushforward measure of $\mathbb P$ through $X$ defined by $\mathbb P_X(B):= \mathbb P(X^{-1}(B))$ for all Borel sets $B \in \mathcal B(\mathbb R)$.
 \end{definition}
 \begin{definition}\label{def:GP}
 A Gaussian process $\big(U(x)\big)_{x\in \mathcal D}$ is a family of Gaussian random variables defined over $(\Omega,\mathcal F, \mathbb P)$ such that for all $n\in\mathbb N$, $(a_1,\dots,a_n) \in \mathbb R^n$ and $(x_1,\dots,x_n)\in \mathcal D$, $\sum_{i=1}^n a_i U(x_i)$ is a Gaussian random variable.
 \end{definition}
 Similar to the Gaussian distribution in finite dimensions, the law the GP induces over the function space $\mathbb R^\mathcal D$ endowed with its product $\sigma$-algebra is uniquely determined by its mean and covariance functions $m(x) = \mathbb E[U(x)]$ and $k(x,x') = \text{Cov}(U(x),U(x'))$. This allows us to write $(U(x))_{x\in \mathcal D} \sim \GP(m,k)$. The covariance function is often called kernel in the literature, in this article the terms are used synonymously. The covariance function $k$ is positive definite over $\mathcal D$, i.e., for all $n\in\mathbb N_{\ge 0}$ and $(x_1,\dots,x_n) \in \mathcal D^n$ the matrix $\big( k(x_i,x_j)\big)_{1\le i,j\le n}$ is positive semidefinite. Conversely, for a given positive definite function over an arbitrary set $\mathcal D$, there exists a centered GP indexed by $\mathcal D$ with this function as its covariance function. For given $\omega \in \Omega$, the sample path $U_\omega \colon \mathcal D \rightarrow \mathbb R$ is the deterministic function $U_\omega(x) := U(x)(\omega)$.

\subsection{Kernel functions}
In the following we recall some kernel functions which are popular in machine learning. For more details and other classes of kernel functions we refer to \cite{Rasmussen2006}.
\paragraph{Squared Exponential} The squared exponential covariance function has the form
$$k_\text{se}(r) = \exp\left(-\frac{r^2}{2\ell^2} \right)$$
with parameter $\ell$ defining the characteristic length-scale. It is infinitely differentiable, which means that the corresponding GP has mean square derivatives of all orders and thus is very smooth. Although, it is the most widely-used kernel within the kernel machine fields, its strong smoothness assumptions are unrealistic for many physical applications.

\paragraph{Mat\'ern class} The Mat\'ern class of covariance functions is given by
\[
k_\text{M}(r) = \frac{2^{1-\nu}}{\Gamma(\nu)} \left( \frac{\sqrt{2\nu} r}{\ell} \right)^\nu K_\nu \left( \frac{\sqrt{2\nu} r}{\ell} \right),
\]
where $\ell$ is again the length scale parameter, $\nu$ allows to controls the smoothness of the function,  $K_\nu$ is the modified Bessel function of the second kind and  $\Gamma(\nu)$ is the Gamma function. Note that for $\nu\rightarrow \infty$ we recover the squared exponential kernel. The process corresponding to the Mat\'ern class is $m$-times mean squared differentiable if and only if $\nu > m$. Note that the kernel becomes especially simple if $\nu = p + 1/2$ for a nonnegative integer $p$. Moreover, setting $\nu = 1/2$ yields the exponential kernel which yields a Gaussian process that is mean squared continuous but not mean squared differentiable. In one dimension, we obtain the Ornstein-Uhlenbeck process with was introduced as a mathematical model of the velocity of a particle undergoing Brownian motion. 
\begin{remark}
For the numerics $k_M(r)$ is multiplied by the signal variance $\sigma^2$. This is an additional parameter allowing to adjust the amplitude of the functions generated by the GP.
\end{remark}

\subsection{Smoothness of sample paths}\label{sec:smoothness}
An overview of recent literature on the smoothness of sample path is given in the introduction of \cite{Henderson2024}. Moreover, a  very general result on the sample path Banach-Sobolev regularity for Gaussian processes can be found in Proposition 3.1 in this article, and a result on Hilbert-Schmidt imbeddings of Sobolev spaces related to GPs is given in its Proposition 4.4. The most relevant result for the applications we have in mind is stated in Example 4.5 of the same article: If $\mathcal D$ is sufficiently smooth, sample paths of a GP corresponding to a Mat\'ern kernel of order $\nu$ lie in $H^m(\mathcal D)$ for $m\in[0,\nu) \cap \mathbb N_0$. Indeed, choosing $\nu = 3/2$ or $\nu=5/2$, which seem to be the most relevant cases in machine learning \cite{Rasmussen2006}, yields sample paths of regularity $H^1(\mathcal D)$ or $H^2(\mathcal D)$, respectively. In addition, it is shown in \cite{Steinwart2019} that in these cases the reproducing kernel Hilbert spaces (RKHS) coincide with $H^1(\mathcal D)$ or $H^2(\mathcal D)$, respectively. However, in general the sample paths of a GP with kernel $k$ do not lie in the RKHS corresponding to the same kernel. For more details and results showing that the sample paths lie in an RKHS corresponding to a power of that kernel in special cases, we refer to \cite{arxivKanagawa}.

\subsection{Marginalization and prediction}\label{sec:marginalization}
The discussion above shows that sampling from a GP prior with appropriate parameters yields functions in Sobolev spaces. To be well-posed, optimization problems in Sobolev spaces usually involve initial value or boundary conditions. We incorporate these with the help of an tailored posterior. Moreover, we exploit the marginalization property which originates from the consistency requirement of GP for the numerical implementation of our scheme. For further details on both, we refer to \cite{Rasmussen2006}. The main point is easy to see in finite dimensions: if the GP specifies 

\[
(y_1,y_2) \sim \mathcal N( \mu,\Sigma) \text{ for } \mu = \begin{pmatrix} a \\ b \end{pmatrix} \in \R^{n_1 + n_2} \text{ and } \Sigma = \begin{pmatrix} A & B \\ B^T & C \end{pmatrix} \in \mathbb R^{(n_1+n_2)\times (n_1+n_2)}
\]
then it holds $y_1 \sim \mathcal N(a, \Sigma_{11})$  and  $y_2 \sim \mathcal N(b, \Sigma_{22})$, respectively. This even holds if $y_2$ contains infinitely many variables, i.e., we can compute the likelihood of $y_1$ under the joint distribution while ignoring $y_2$. This helps us particularly in the implementation, since we can compute finite dimensional projections of the underlying infinite object on our computer. Indeed, we will represent the functions (sample paths of the GP) by their finite dimensional projections which will correspond to the evaluations on a predefined mesh.

\begin{figure}
    \centering
    \includegraphics[width=0.95\linewidth,height=0.4\linewidth]{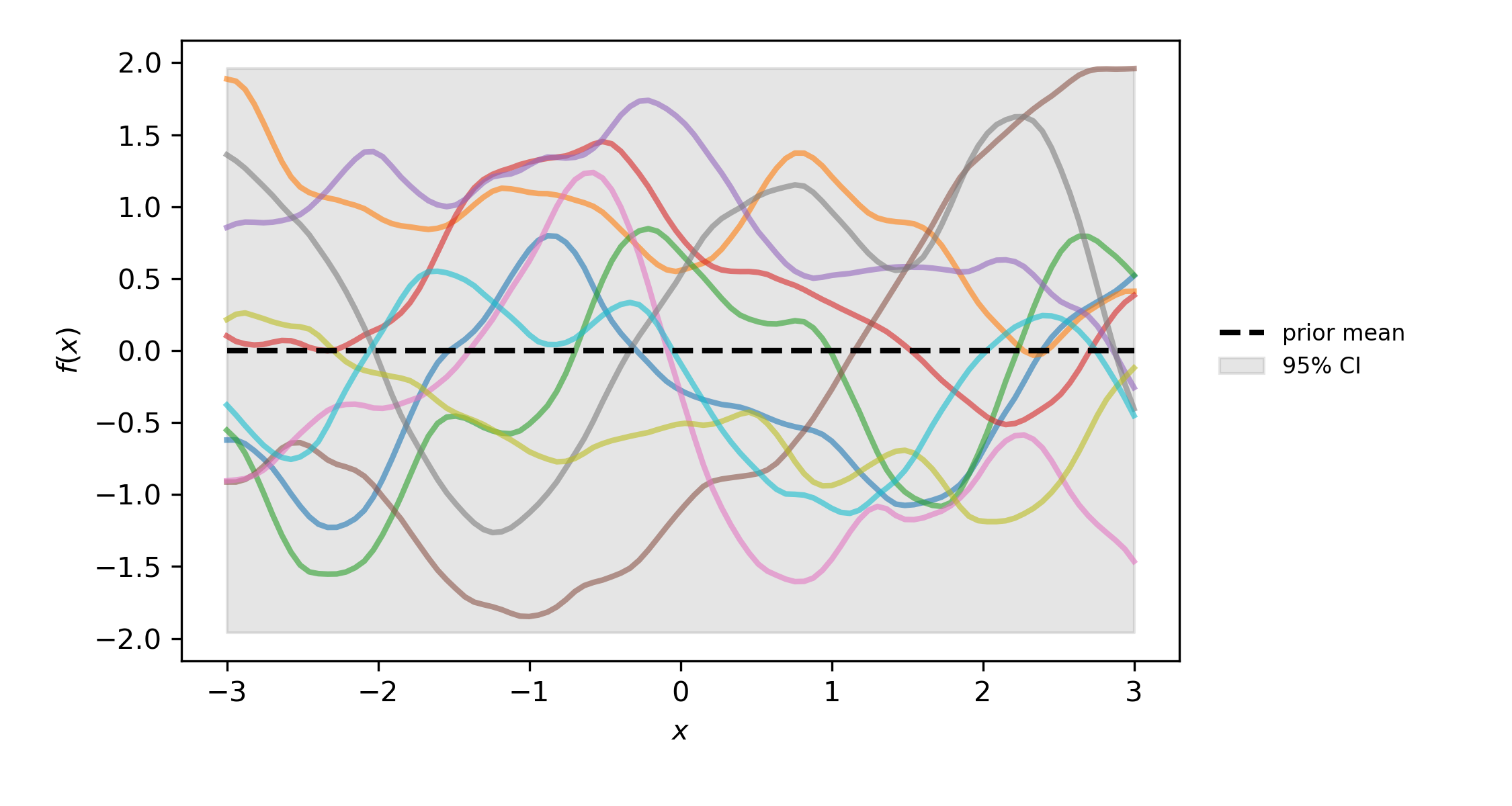}
    \includegraphics[width=0.95\linewidth,height=0.4\linewidth]{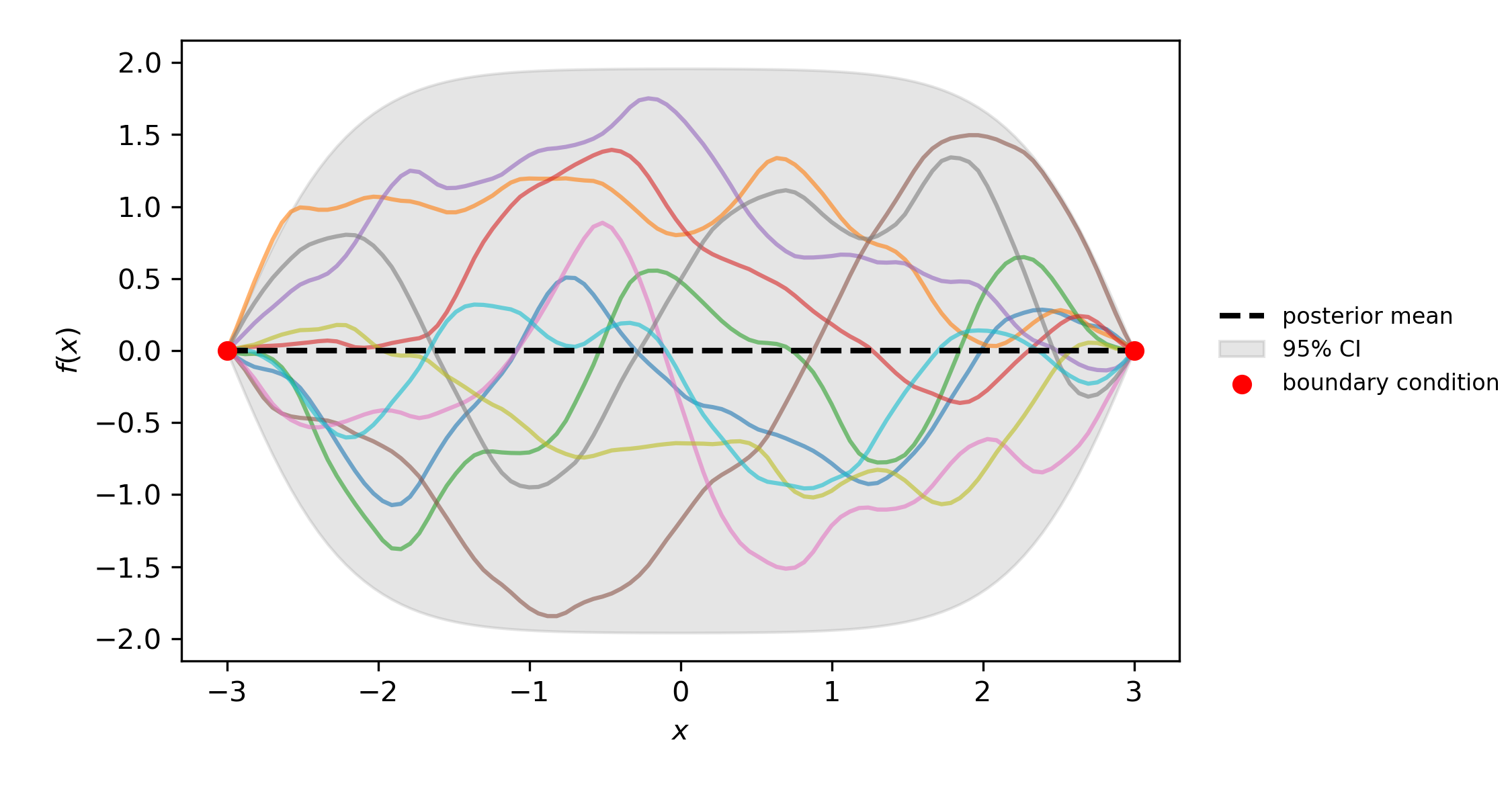}
    \includegraphics[width=0.95\linewidth,height=0.4\linewidth]{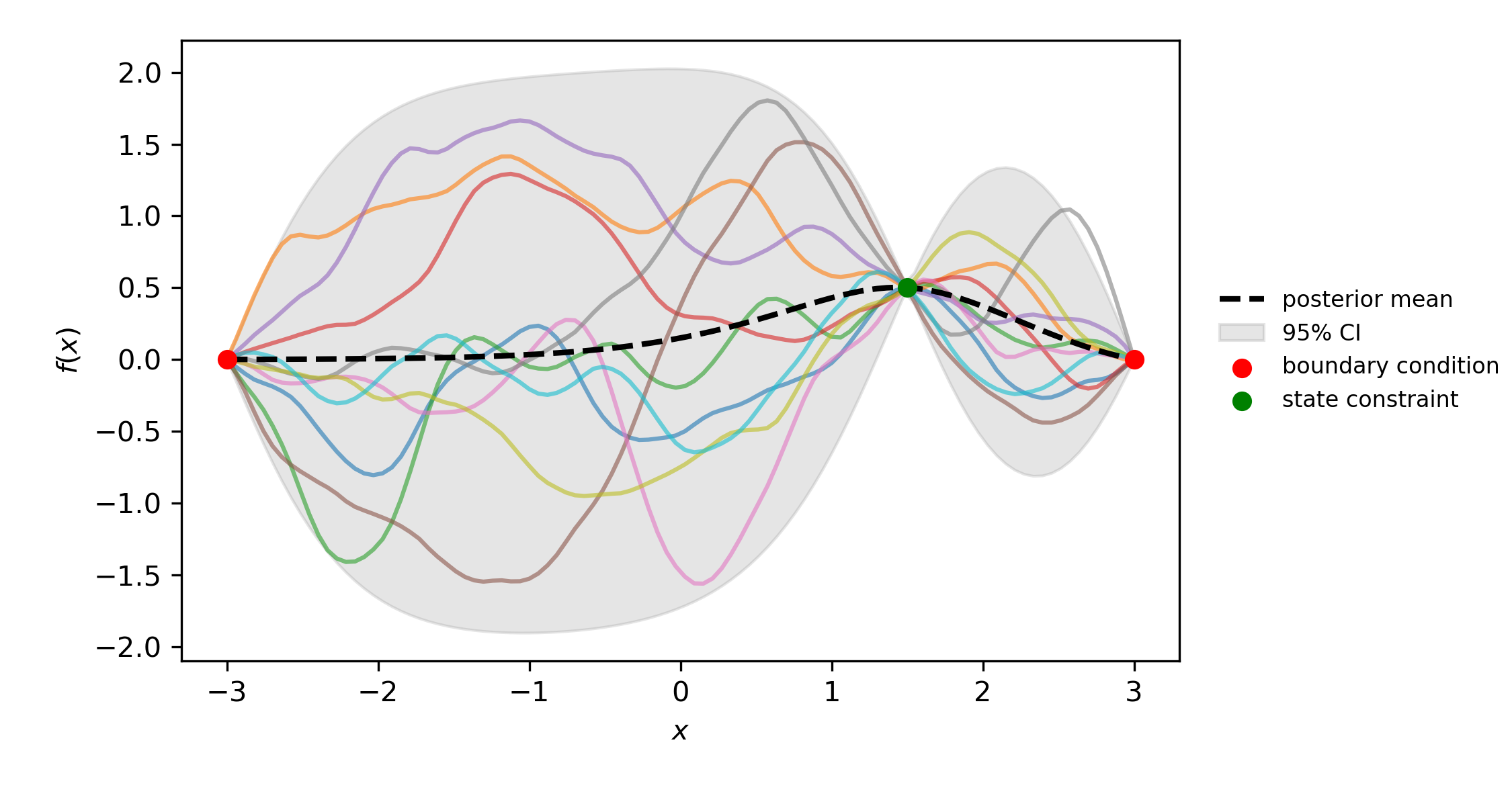}
    \caption{Visualization of samples drawn from the GP defined via Mat\'ern Kernel with parameters $\nu=2.5$, $\ell=1$ and signal variance $\sigma^2=1$.  The shaded areas indicate the 95\% confidence band. Top: $10$ samples drawn from prior. Middle: $10$ samples drawn from posterior which takes the boundary conditions at $x=-3$ and $x=3$ into account. Bottom: Visualization of $10$ samples drawn from the GP with additional state constraints.
    }
    \label{fig:prior_posterior_state}
\end{figure}

Similarly, we will handle given initial or boundary data, and if applicable also state constraints,  specified by the application at hand. In fact, we split the variables into thee parts $(y_1,y_2,y_3)$, where $y_1$ are the projections as above, $y_2$ contains information given at data points and $y_3$ is infinite dimensional. As above $y_3$ can be ignored in the computation. In order to incorporate the data of $y_2$ in the finite projection $y_1$, we employ Bayes rule as it is usually done when GP are used for predictions
\[ p(y_1 | y_2) = \frac{p(y_1,y_2)}{p(y_2)}. \]
Since we are in the Gaussian setting, we can compute explicitly
\begin{equation}\label{eq:posterior}
p(y_1 | y_2) = \mathcal N(a + BC^{-1}(y_2-b), A-BC^{-1}B^\top), 
\end{equation}
where we used the notation from above. It is common to choose mean zero, that is, $a$ and $b$ zero vectors \cite{Rasmussen2006} and we will do so in the following. A visualization of samples from the prior, samples from the posterior with boundary conditions and samples from the posterior with additional state constraints are shown in Figure~\ref{fig:prior_posterior_state}. The shaded area shows the confidence band which is computed pointwise by adding the mean and two time the standard deviation at every position.

\section{Algorithm and implementation}

We note that the dynamics proposed in \eqref{eq:timediscrete} is discrete in time but continuous in space. Indeed, $X_j^i$ is a function from some Sobolev space for each $j$ and each $i$. To obtain a scheme that is feasible for implementation, we exploit the well-known marginalization property, see Section~\ref{sec:marginalization}.

\subsection{Algorithm in function space}
Let us denote by $\GP(m,k,\theta)$ the Gaussian Process with mean $m$, kernel $k$ and additional parameters $\theta$, i.e.,~$\theta = (\ell,\nu)$ in case of the Mat\'ern class. Moreover, let us encode initial value, boundary conditions or additional constraints with the help of an operator $c(u)=0$. 
For notational convenience, we denote by $\GP_c(m,k,\theta)$ the posterior set. In more detail, $U \in \GP_c(m,k,\theta)$ is from the same space as $\GP(m,k,\theta)$ with the additional property $c(U)=0$, that means, that initial value conditions, boundary conditions and additional constraints are satisfied. For notational convenience, we define the set $\GP_0(n,k,\theta)$ as posterior set, where all constraints are assumed to be homogeneous. For example, if $U \sim \GP_c(m,k,\theta)$ satisfies an initial condition $U(t_0) = \bar u \ne 0$, a sample $\xi\sim \GP_0(m,k,\theta)$ would satisfy $\xi(t_0)=0$. 
In every iteration of the algorithm we consider the update
\begin{align*}
&U^i\gets U^i-\lambda (U^i - v_f^\alpha(U))\tau + \sqrt{2\tau} \|U^i - v_f^\alpha(U)\|_{\mathcal U} \,\xi^i, \\
&\qquad \quad U_0 \in \GP_c(m,k,\theta), \qquad \xi^i \in \GP_0(m,k,\theta).
\end{align*}
Since $v_f^\alpha$ is a weighted mean, the difference $U^i - v_f^\alpha(U)$ satisfies homogeneous conditions. By definition, $\xi^i$ satisfies homogeneous conditions. Altogether, this ensures that the initial, boundary and additional conditions of the respective problem are met by each iterate of the algorithm, moreover, that the regularity of the functions is preserved in every iteration. For more details on the update procedure we refer to Algorithm~\ref{alg:GP-CBO}.

\begin{algorithm}
\caption{GP-CBO}
\label{alg:GP-CBO}
\begin{algorithmic}[1] 
\Require number of CBO agents $N$; constraint posterior $\GP_c(m,k,\theta)$; homogeneous posterior $\GP_0(m,k,\theta)$; time horizon $T$, time step size $\tau$; weight parameter $\alpha$, drift parameter $\lambda$; 
\State $J \gets int(T/\tau)$
\State $U_0=(U^i)_i  \gets $ $N$ independent samples $U^i \sim \GP_c(m,k,\theta)$

\For{$j=0,\dots J-1$}
    \State evaluate $v_f^\alpha(U_j)$ 
    \State $\xi=(\xi^i)_i \gets $ $N$ independent samples $\xi^i \sim \GP_0(m,k,\theta)$%, $i=1,\dots,N$
    \State compute the weighed mean $v_f^\alpha(U_j)$
    \For{$i=1,\dots N$}
    \State $U_{j+1}^i \gets U_{j}^i-\lambda (U_j^i - v_f^\alpha(U_j))\tau + \sqrt{2\tau} \|U_j^i - v_f^\alpha(U_j)\|_{\mathcal U} \,\xi_j^i$
    \EndFor
\EndFor
\State \Return $v_f^\alpha(U_{J})$, $U_{J}$
\end{algorithmic}
\end{algorithm}
\noindent We want to emphasize that, by construction of $\GP_c$, given data, that is, state constraints, initial conditions or boundary conditions are included in $\mathcal U$. If we consider, for example, an initial value problem given by 
\begin{equation*}
\frac{d}{dt} u = g(u), \qquad u(t_0) = u_0.
\end{equation*}
We could apply Algorithm~\ref{alg:GP-CBO} with cost functional $f(u) = \|\frac{d}{dt} u - g(u) \|_{L^2(0,T)}^2$ and admissible set
$\mathcal U = \{ u\in H^1(0,T) \colon u(t_0) = u_0 \}$. Choosing appropriate parameters for the GP and sampling from the posterior $\GP_c$ then ensures that the samples $U_0$ satisfy $U_0 \in \mathcal U$.

\begin{remark}\emph{(Isotropic versus anisotropic noise)} Recent articles on CBO variants in finite dimensions discuss two different types of noise terms: isotropic given by $|x-v_f|dB_t^i$ and anisotropic $\emph{diag}(x-v_f) dB_t^i$ as first proposed in \cite{carrillo2021anisotropic}. In the function space setting, the isotropic case translates to the noise introduced above, and the anisotropic noise corresponds to pointwise multiplication of the functions $(U^i_j - v_f^\alpha(U_j)) \xi_j^i$. Noting that the pointwise multiplication may lead to loss of regularity (for example for $f,g\in L^2$, in general $fg \notin L^2$), we chose the isotropic case. This choice is further underpinned by numerical simulations we ran, which showed slower or even premature convergence for the algorithm with pointwise multiplication in the noise term.
\end{remark}

\subsection{Discretization and implementation}
Each $U_j$ in Algorithm~\ref{alg:GP-CBO} is a $N$-dimensional vector of samples from $\GP_c(m,k,\theta)$, in particular, each $U_j^i$ is a function. To obtain an implementable version, we need to discretize these functions. Let therefore $\mathcal X = (x_p)_{p\in \mathcal P}$ be a discretization of the domain of $U_j^i$, where $\mathcal P$ denotes the corresponding index set that depends on the spatial dimension of the domain. Evaluating the mean $m$ and the kernel $k$ along $\mathcal X$, we obtain tensors $m_{|\mathcal X}$ and $k_{|\mathcal X}$. 

In the particular case where the domain of the function is an interval discretized with $P\in\mathbb N$ points, we obtain $\mathcal P=\{1,\dots,P\}$, which yields 
\[
m_{|\mathcal X} = (m(x_p))_p \in \mathbb R^P, \qquad k_{|\mathcal X} = (k(x_p,x_q))_{p,q} \in \mathbb R^{P\times P}.
\]

Let $V \sim \GP(m,k)$ a sample path of the prior. A discretization of $V$ is a multivariate Gaussian random variable
\[
v \sim \mathcal N(m_{|\mathcal X}, k_{|\mathcal X} + \sigma_\text{GP}^2 I),
\]
where $I$ denotes the identity matrix and $\sigma_\text{GP}^2$ is a noise parameter of the Gaussian process. Figure~\ref{fig:prior_posterior_state} (top) illustrates $10$ samples $v$ from the prior.

To include initial values, boundary conditions or state state constraints, we compute the corresponding posterior distribution with the help of Bayes rule. Since the initial data, boundary data and state constraints are known information, these play the role of training points in the usual setting of GPs. To discuss the details, we denote by $x_\text{in}$ the points where initial data are prescribed, $x_\text{bc}$ the points with prescribed boundary conditions and $x_\text{sc}$ the points at which we have state constraints. For notational convenience we discuss the case where all types of conditions are present. If one case does not apply, we assume the corresponding vectors to be empty. The vector $x_\text{train} = (x_\text{in},x_\text{bc},x_\text{sc}) \in \R^{d_\text{data}}$ collects all the given data in a vector of the dimension $d_\text{data}$.

Now, we build the matrices needed for the computation of the posterior as follows
\begin{align*}
A:= k_{|\mathcal X} + \sigma^2_\text{GP}I \in \R^{P\times P}, \quad B = k_{|\mathcal X,x_\text{train}} \in \R^{P\times d_\text{data}}, \quad C := k_{|x_\text{train}} \in \R^{d_\text{data}\times d_\text{data}},
\end{align*}
by pairwise evaluation of the kernel with the respective points. Note that $A$ corresponds to the covariance matrix we already constructed for the prior. In $B$ we assemble the values combining the information of mesh and training points. And $C$ collects the evaluations of the training points. The distribution of the posterior is then given by \eqref{eq:posterior}, hence in line $2$ of Algorithm~\ref{alg:GP-CBO} we draw random variables 
\[
U^i \sim GP_c(m,k,\theta) = \mathcal N (BC^{-1}x_\text{train}, A-BC^{-1}B^\top).
\]
To preserve the prescribed data in all iterations, we make sure that the updates in line $5$ of Algorithm~\ref{alg:GP-CBO} are set to zero at all training points, in more details the random variables satisfy
\[
\xi^i \sim GP_0(m,k,\theta) = \mathcal N(0, A-BC^{-1}B^\top).
\]

\section{Applications}\label{sec:problemClasses}
The algorithm specified in the previous section can be used in very general settings. In the following we illustrate its performance with the help of boundary value problems with state constraints and nonlinear optimal control problems constrained by ordinary differential equations.

\subsection{Boundary value problems (BVP)} 
Since unconstrained BVP are well-understood and in special cases even admit analytic solutions, we choose them for the first tests of the proposed algorithm. Indeed, we consider problems of the form
$e(u)=0$, where the operator $e\colon \mathcal U \rightarrow Z$ may involve  derivatives of $u$ and optional state constraints. To guarantee well-posedness, boundary values have to be specified. Most common are boundary conditions of Dirichlet, Neumann or Robin-type, we use Dirichlet conditions in the following as they are the simplest case to realize with the help of GP. As the constraints are incorporated in the posterior distribution, we are free to choose an appropriate cost functional without explicitly addressing the constraints. Therefore a candidate for an appropriate cost functional is given by
\[
f(u) = \| e(u) \|_Z^2 
\]
For example, $Z=L^2(\mathcal D)$ yields $f(u) = \int_{\mathcal D} e(u)^2 dx $.

\begin{remark}
For BVP which allow for the well-posedness analysis based on energy estimates, the energy functional may be used as alternative cost functional. Note that these usually require less regularity and hence less derivatives of $u$ have to be computed to evaluate the cost functional. 
\end{remark}

\subsection{Optimization problems constrained by differential equations}
The second class of problems we consider are global optimal control problems. We think of constraints given by ordinary differential equations (ODE) or partial differential equations (PDE) which often arise in physics. In more detail, we consider problems of the form
\begin{equation}\label{prob:OC}\tag{OC}
\min\limits_{y\in \mathcal Y, u\in \mathcal U} \mathcal J(y,u) \quad \text{subject to} \quad e(y,u)=0,
\end{equation}
where $\mathcal J \colon \mathcal Y \times \mathcal U \rightarrow \mathbb R$ is a cost functional and the operator $e\colon \mathcal Y \times \mathcal U \rightarrow Z$ encodes the state constraint. In more detail $e(y,u)=0$ implies that $y$ is the solution of the differential equation for given control $u$. Assuming the well-posedness of the state equation, it is common practice to introduce the control-to-state map $\mathcal S \colon \mathcal U \rightarrow \mathcal Y$. If $y=\mathcal S(u)$, we introduce the reduced cost functional 
\[
f(u) := \mathcal J(\mathcal S(u),u).
\]
Then \eqref{prob:OC} problem can be rewritten as 
\[
\min\limits_{u \in \mathcal U} f(u). 
\]
Hence, instead of considering the constraint explicitly, it is taken into account via $\mathcal S$. However, we are now in an unconstrained setting and may apply Algorithm~\ref{alg:GP-CBO}.

\begin{remark}
Existence results for \eqref{prob:OC} are usually based on coercivity of $\mathcal J$ w.r.t.~$u$, weak continuity and boundedness of the map $\mathcal S$ and weakly lower semicontinuity of $\mathcal J$. The arguments are therefore not constructive and often first-order optimality conditions are used to characterize and compute local minimizers. We refer to \cite{Hinze2009} for more details. In contrast, Algorithm~\ref{alg:GP-CBO} is a zero-order method that aims to approximate the global solution of \eqref{prob:OC}. 
This is particularly of interest if the constraints or the cost functional are nonlinear, since this causes usually nonconvexity and therefore local and global solutions may differ.
\end{remark}

\section{Numerical results}
In the following we show numerical results for the different problem classes discussed in Section~\ref{sec:problemClasses}. We begin with boundary value problems. For the numerical computations we used the Mat\'ern kernel with length scale $\ell = 1$ and smoothness parameter $\nu = 2.5$ if not stated otherwise. The code used for the numerical results is available here \cite{zenodo_folder}.

\subsection{Boundary value problems}
First, we test our approach in one and two dimensional settings with  boundary value problems that admit exact solutions. Then we add some state constraints, which pose a mayor challenge for standard solution methods and justify the additional computational effort of the CBO approach. 
\subsubsection{One dimensional case}
Indeed, we consider the following boundary value problem \cite{Simmons2017}:

\begin{equation}
    u''(x) + u(x) = 0
\end{equation}
with the boundary conditions
\begin{equation}
    u(0) = 0, \quad u\left(\frac{\pi}{2}\right) = 2.
\end{equation}

In the CBO approach we measure feasibility of a candidate function $u$ via the residual of the differential equation and use the cost functional
\begin{equation}
    f(u) \;=\; \int_{0}^{\pi/2} \bigl(u''(x)+u(x)\bigr)^2\,dx,
\end{equation}
where the boundary conditions are enforced by sampling from a Gaussian process posterior conditioned on $u(0)=0$ and $u(\pi/2)=2$.

\begin{figure}[ht]
    \centering
    \begin{minipage}{0.48\textwidth}
        \centering
        \includegraphics[width=\textwidth]{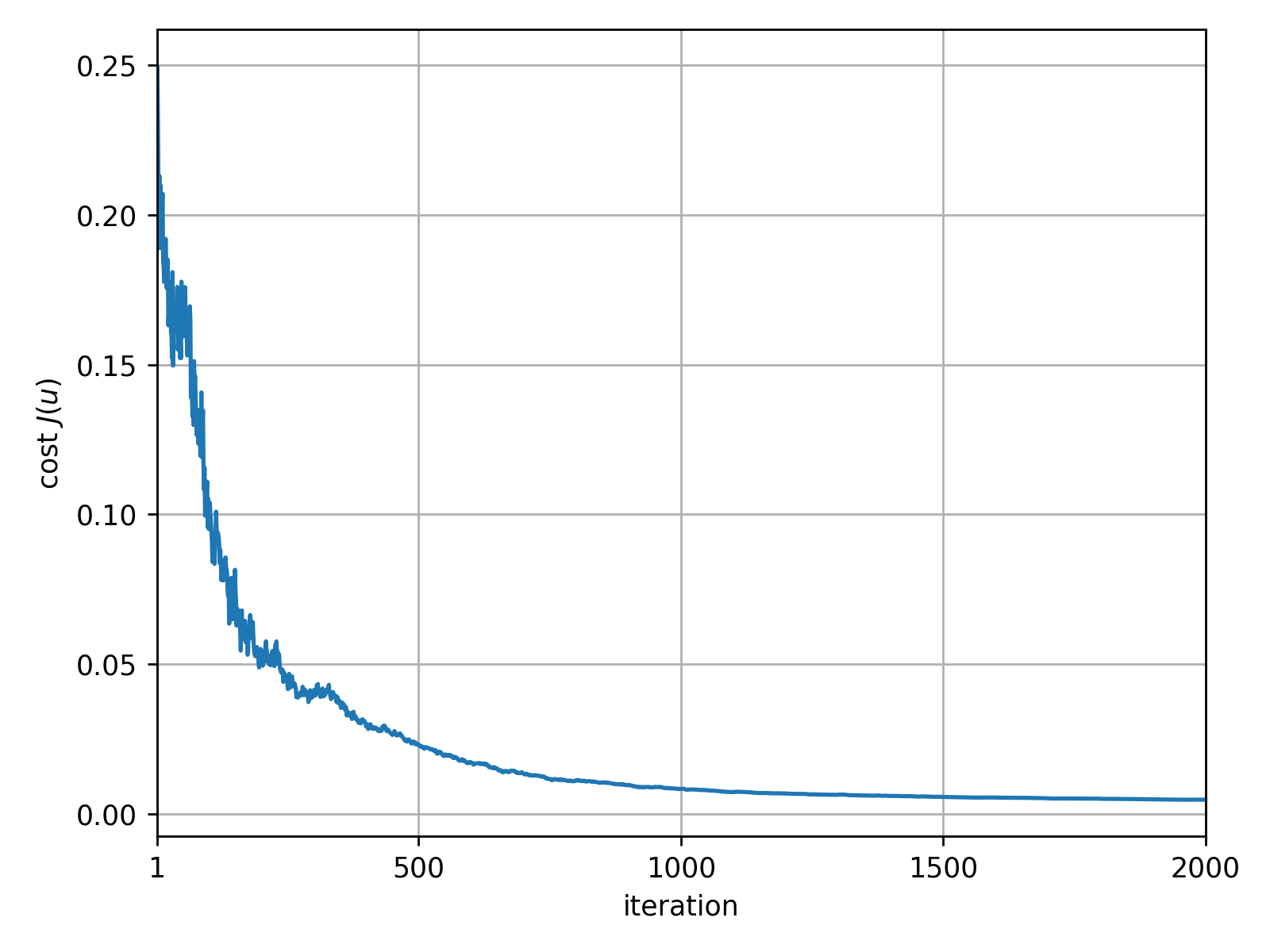}
        \caption*{(a) Cost functional $J(u)$}
    \end{minipage}
    \hfill
    \begin{minipage}{0.48\textwidth}
        \centering
        \includegraphics[width=\textwidth]{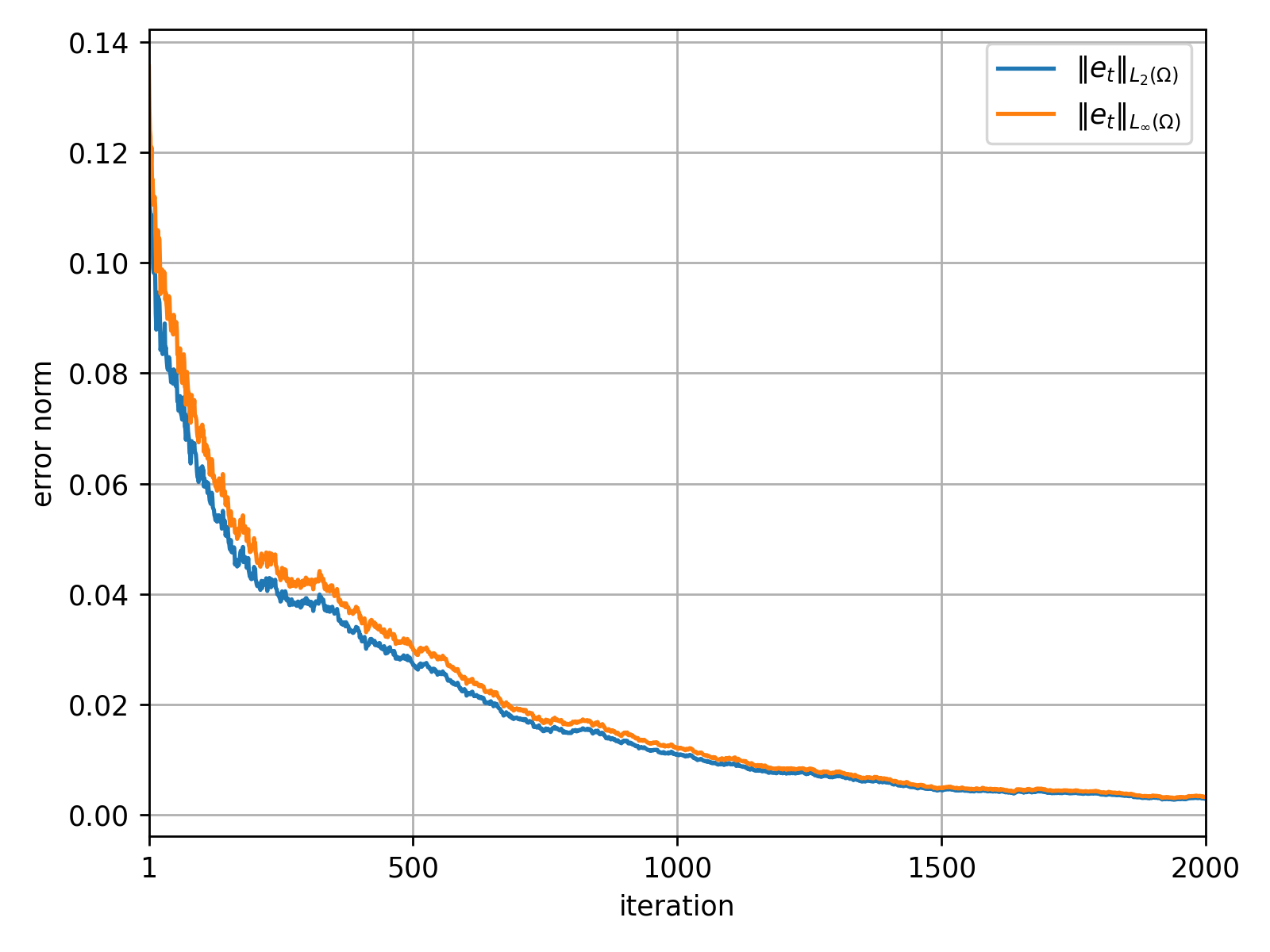}
        \caption*{(b) Error norms $\|e\|_{L_2(\Omega)}$ and $\|e\|_{L_\infty(\Omega)}$}
    \end{minipage}
    \caption{Evolution of the cost functional and the error norms during the CBO iterations.}
    \label{fig:cbo_convergence}
\end{figure}

The exact solution is given by
\begin{equation}
    u_{*} = 2\sin(x).
\end{equation}
To quantify the difference between the CBO solution $u$ and the exact solution $u_{*}$,
we consider the error function
\begin{equation}
e(x) = u(x) - u_{*}(x).
\end{equation}
We measure this error in the $L_2$-norm and $L_\infty$-norm, respectively defined by
\begin{equation}
\|e\|_{L_2(\Omega)} := \left( \int_{\Omega} |e(x)|^2 \, dx \right)^{1/2},
\qquad
\|e\|_{L_\infty(\Omega)} := \sup_{x \in \Omega} |e(x)|,
\end{equation}
where $\Omega = [0,\pi/2]$.
During the CBO iterations, we monitor the evolution of the cost functional and the error norms.
Figure~\ref{fig:cbo_convergence} summarizes the convergence behavior: the left panel shows the
decay of the cost functional $J(u)$ over the iterations, while the right panel presents the
evolution of the $L_2$- and $L_\infty$-norms of the error $e(x)$.

\begin{figure}
    \centering
    \includegraphics[width=0.45\linewidth]{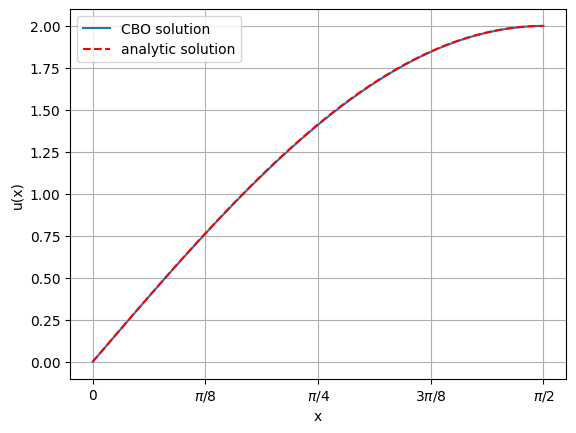}
    \includegraphics[width=0.45\linewidth]{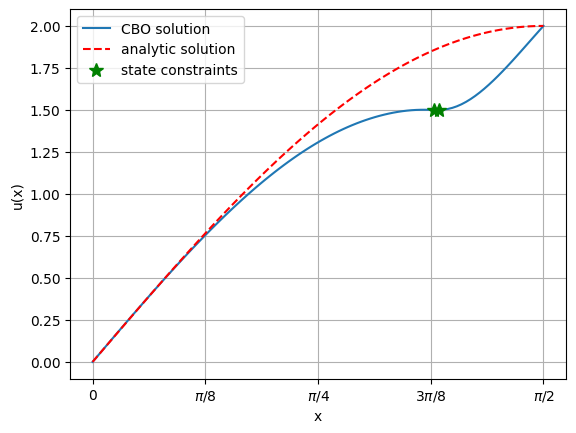}
    \captionof{figure}{Comparison between unconstrained and constrained solution found by CBO. }
    \label{fig:BVP1d_comparison}
\end{figure}

Now, we consider the problem with additional state constraints at two points on the mesh
\[
x_1 = 1.189997 \quad \text{with} \quad y(x_1) =1.85673, \qquad
x_2 =1.20586 \quad \text{with} \quad y(x_2) = 1.86829.
\]
Due to the additional constraints the problem, we do not have analytical solution, however we see in Figure~\ref{fig:BVP1d_comparison} that the solution found by CBO behaves similar to the unconstrained case.

\subsubsection{Two dimensional case}
We extend our tests to 2d to solve the following Poisson problem:
\begin{subequations}\label{eq:poisson2d}
\begin{align}
- \Delta u(x,y) &= -6, \quad (x,y) \in \Omega, \\ 
u(x,y) &= 1 + x^2 + 2y^2, \quad (x,y) \in \partial \Omega, \\ 
\Omega &= [0,1] \times [0,1]. 
\end{align}
\end{subequations}

The exact solution of problem \eqref{eq:poisson2d} is
\begin{equation}
u_e(x,y) = 1 + x^2 + 2y^2. \label{eq:poisson2d_exact}
\end{equation}
The computational domain $\Omega = [0,1] \times [0,1]$ is discretized by a uniform structured grid with 
$N_x = N_y = 30$ nodes in each spatial direction.

\begin{figure}[H]
    \centering
    \includegraphics[width=0.7\textwidth]{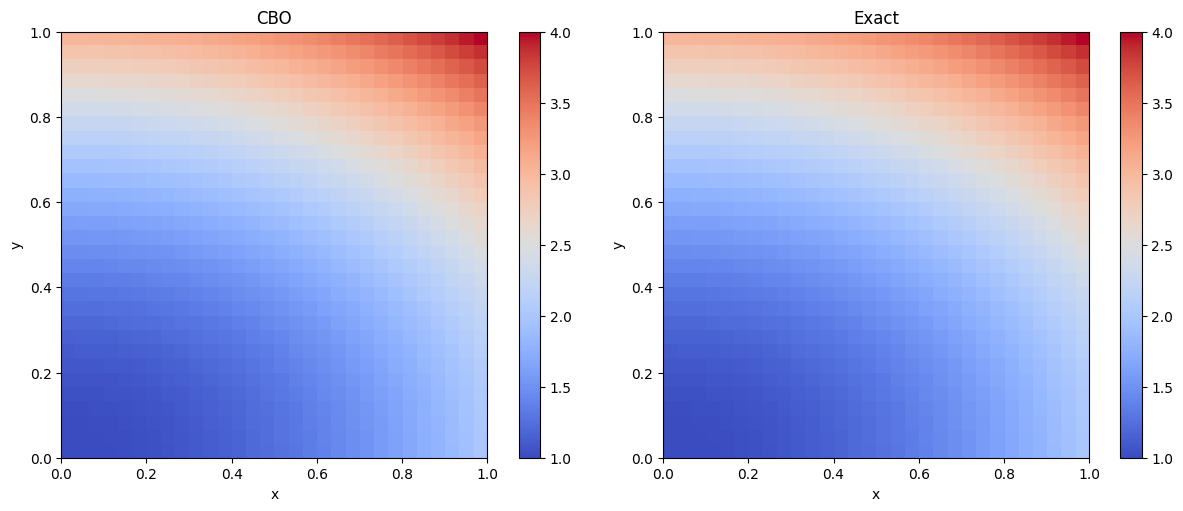}
    \caption{Comparison between the CBO solution and the analytical solution for the two-dimensional Poisson problem \eqref{eq:poisson2d}.}
    \label{fig:cbo_vs_exact_2d_unconstrained}
\end{figure}

The numerical solution obtained by CBO is compared with the analytical solution in Figure~\ref{fig:cbo_vs_exact_2d_unconstrained}.
In the CBO approach, feasibility of a candidate function $u$ is measured via the residual of the Poisson equation, and we use the cost functional
\begin{equation}
    J(u) \;=\; \int_{\Omega} \bigl(\Delta u(x,y) - 6\bigr)^2 \, dx\,dy,
\end{equation}
where $\Delta u = u_{xx}+u_{yy}$ and $\Omega=[0,1]\times[0,1]$.

To quantify the difference between the CBO solution $u$ and the exact solution $u_e$,
we define the error function
\begin{equation}
e(x,y) = u(x,y) - u_e(x,y),
\end{equation}
and measure it in the $L_2$-norm and $L_\infty$-norm, respectively,
\begin{equation}
\|e\|_{L_2(\Omega)} := \left( \int_{\Omega} |e(x,y)|^2 \, dx\,dy \right)^{1/2},
\qquad
\|e\|_{L_\infty(\Omega)} := \sup_{(x,y)\in \Omega} |e(x,y)|.
\end{equation}

\begin{figure}[H]
    \centering
    \begin{minipage}{0.48\textwidth}
        \centering
        \includegraphics[width=\textwidth]{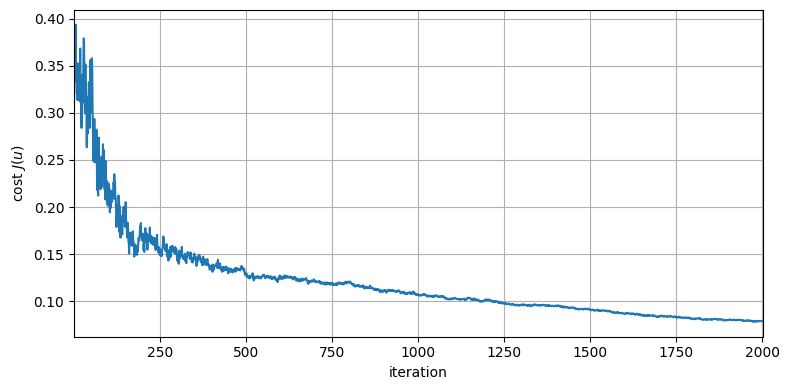}
        \caption*{(a) Cost functional $J(u)$}
    \end{minipage}
    \hfill
    \begin{minipage}{0.48\textwidth}
        \centering
        \includegraphics[width=\textwidth]{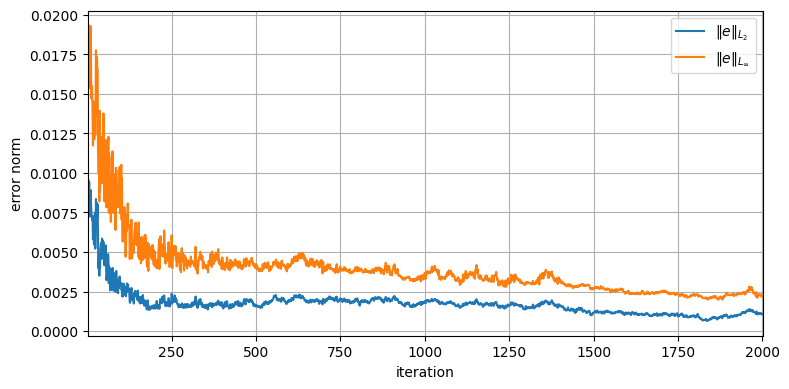}
        \caption*{(b) Error norms $\|e\|_{L_2(\Omega)}$ and $\|e\|_{L_\infty(\Omega)}$}
    \end{minipage}
    \caption{Evolution of the cost functional and the error norms during the CBO iterations for the two-dimensional Poisson problem \eqref{eq:poisson2d}.}
    \label{fig:cbo_convergence_2d}
\end{figure}

\begin{figure}[H]
    \centering
    \begin{minipage}{0.48\textwidth}
        \centering
        \includegraphics[width=\textwidth]{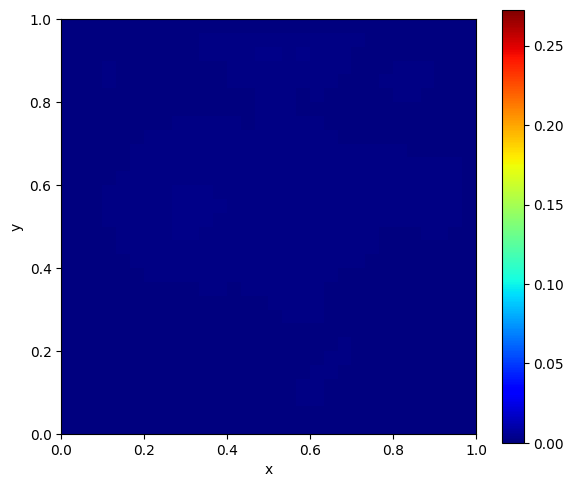}
        \caption*{(a) Unconstrained case: heat map of $\lvert u_{\mathrm{CBO}} - u_{\mathrm{exact}} \rvert$.}
    \end{minipage}
    \hfill
    \begin{minipage}{0.48\textwidth}
        \centering
        \includegraphics[width=\textwidth]{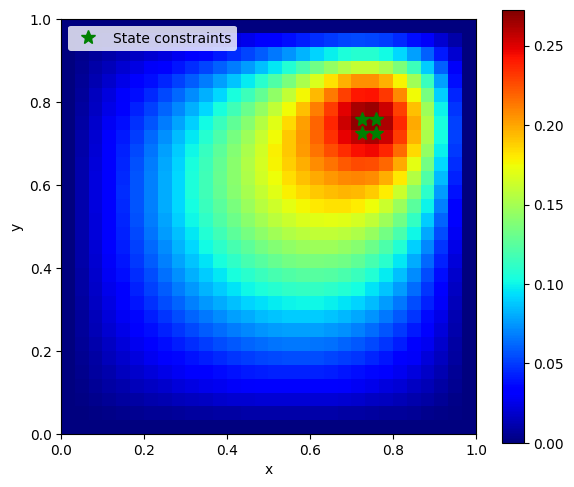}
        \caption*{(b) Constrained case: heat map of $\lvert u_{\mathrm{CBO}} - u_{\mathrm{exact}} \rvert$.}
    \end{minipage}
    \caption{Comparison of the absolute error $\lvert u_{\mathrm{CBO}} - u_{\mathrm{exact}} \rvert$
    between the unconstrained and constrained CBO solutions for the two-dimensional Poisson problem \eqref{eq:poisson2d}.}
    \label{fig:absdiff_2d_comparison}
\end{figure}

During the CBO iterations, we monitor the evolution of the cost functional and the error norms, shown in
Figure~\ref{fig:cbo_convergence_2d}: the left panel displays the cost $J(u)$ over the iterations, while the right panel reports the evolution of $\|e\|_{L_2(\Omega)}$ and $\|e\|_{L_\infty(\Omega)}$.
Now, we consider the problem with additional state constraints at four points on the mesh
\[
\begin{aligned}
(x_1,y_1) &= (0.72413793,\,0.72413793), \\
(x_2,y_2) &= (0.75862069,\,0.72413793), \\
(x_3,y_3) &= (0.72413793,\,0.75862069), \\
(x_4,y_4) &= (0.75862069,\,0.75862069),
\end{aligned}
\]
with prescribed values
\[
\begin{aligned}
u(x_1,y_1) &= 2.31581451, \\
u(x_2,y_2) &= 2.36183115, \\
u(x_3,y_3) &= 2.40784780, \\
u(x_4,y_4) &= 2.45386445.
\end{aligned}
\]
Due to the additional state constraints, we do not have an analytical solution for the constrained
problem. Instead, we compare the
absolute difference
\[
\text{abs\_diff}(x,y) = \bigl|u_{\mathrm{CBO}}(x,y) - u_{\mathrm{exact}}(x,y)\bigr|
\]
for the unconstrained and constrained cases, depicted in Figure~\ref{fig:absdiff_2d_comparison}
in the form of heat maps.

\begin{figure}[H]
    \centering
    \includegraphics[width=0.4\textwidth]{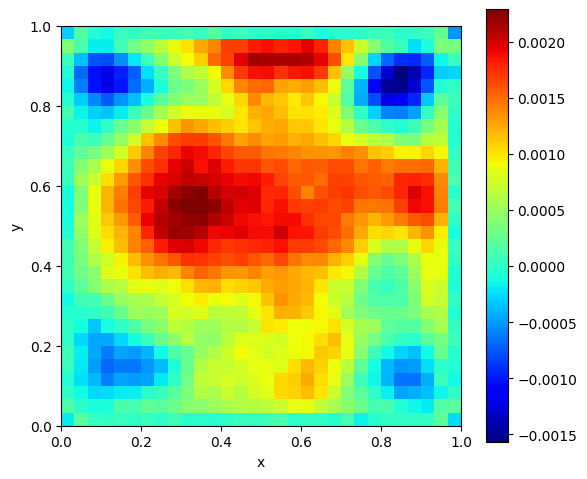}
    \caption{Signed difference $\text{diff}(x,y) = u_{\mathrm{CBO}}(x,y) - u_{\mathrm{exact}}(x,y)$
    for the unconstrained Poisson problem \eqref{eq:poisson2d}.}
    \label{fig:tiny_difference_2d}
\end{figure}

 We use the same color bar range in Figures~\ref{fig:absdiff_2d_comparison}(a) and \ref{fig:absdiff_2d_comparison}(b) in order to ensure a fair comparison between the two configurations, which does not mean that the heat map in Figure~\ref{fig:absdiff_2d_comparison}(a) is identically zero everywhere. To further illustrate that the unconstrained CBO solution is not exactly equal to the analytical one, we also plot the signed difference
\[
\text{diff}(x,y) = u_{\mathrm{CBO}}(x,y) - u_{\mathrm{exact}}(x,y)
\]
for the unconstrained case; see Figure~\ref{fig:tiny_difference_2d}. This plot reveals small but nonzero deviations across the domain, which are too small to be clearly visible when using the shared color scale in Figure~\ref{fig:absdiff_2d_comparison}.

We study the performance of CBO on a more complex problem, namely a two-dimensional nonlinear example. In order to obtain an analytical solution for comparison, we make use of the \emph{method of manufactured solutions}, as described in \cite{fenics}. Hence we adapt problem \eqref{eq:poisson2d} to the nonlinear equation
\begin{equation}
- \Delta u + u^3 = f(x,y), 
\quad (x,y) \in \Omega,
\label{eq:nonlinear_pde}
\end{equation}
where the source term is chosen as
\begin{equation}
f(x,y) = -6 + (1 + x^2 + 2y^2)^3
\label{eq:nonlinear_source}
\end{equation}
and 
Dirichlet boundary condition  as in 
\eqref{eq:poisson2d}, namely
\[
u(x,y) = 1 + x^2 + 2y^2 
\quad \text{on } \partial \Omega.
\]

\begin{figure}[H]
    \centering
    \includegraphics[width=0.8\textwidth]{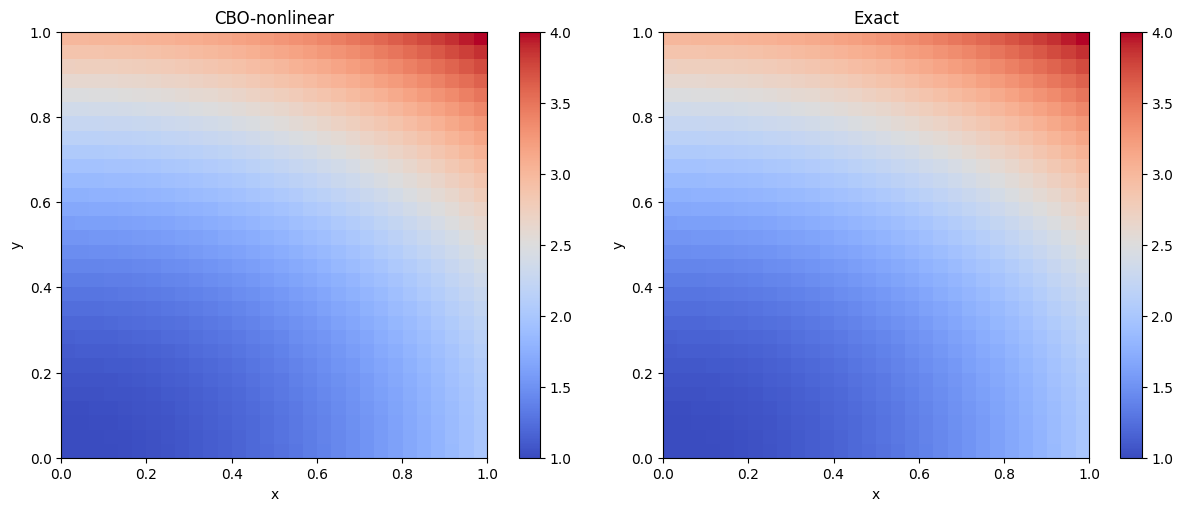}
    \caption{Comparison between the CBO solution for \eqref{eq:nonlinear_pde} and the exact solution.}
    \label{fig:nonlinear_comparison}
\end{figure}

\begin{figure}[H]
    \centering
    \begin{minipage}{0.48\textwidth}
        \centering
        \includegraphics[width=\textwidth]{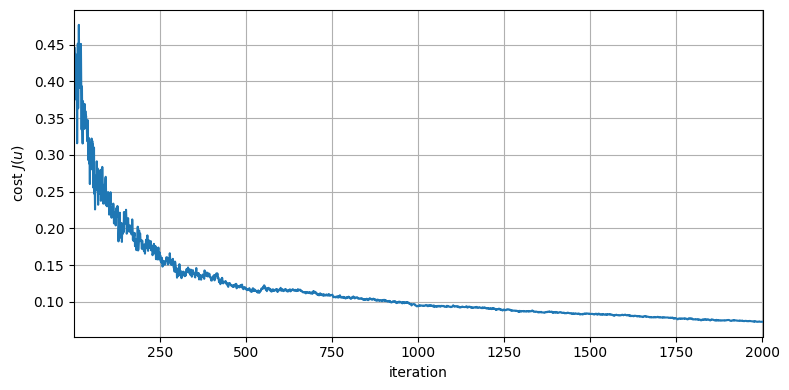}
        \caption*{(a) Cost functional $J(u)$}
    \end{minipage}
    \hfill
    \begin{minipage}{0.48\textwidth}
        \centering
        \includegraphics[width=\textwidth]{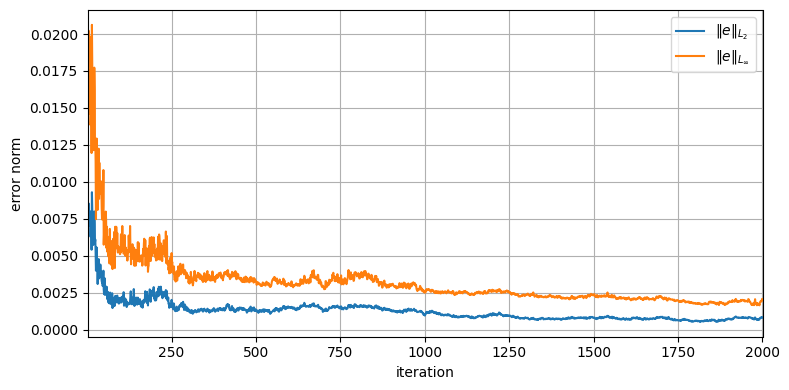}
        \caption*{(b) Error norms $\|e\|_{L_2(\Omega)}$ and $\|e\|_{L_\infty(\Omega)}$}
    \end{minipage}
    \caption{Evolution of the cost functional and the error norms during the CBO iterations for the two-dimensional nonlinear problem \eqref{eq:nonlinear_pde}.}
    \label{fig:cbo_convergence_nonlinear}
\end{figure}

By construction, we obtain the same exact solution as in 
\eqref{eq:poisson2d_exact},
\[
u_e(x,y) = 1 + x^2 + 2y^2.
\]
The output of the nonlinear CBO solution compared to the exact solution is shown in 
Figure~\ref{fig:nonlinear_comparison}. 
During the CBO iterations, we monitor the evolution of the cost functional and the error norms, shown in
Figure~\ref{fig:cbo_convergence_nonlinear}: the left panel displays the cost $J(u)$ over the iterations, 
while the right panel reports the evolution of $\|e\|_{L_2(\Omega)}$ and $\|e\|_{L_\infty(\Omega)}$.
We can see that the performance of CBO was not influenced by the increased complexity of the nonlinear example, 
since for CBO this corresponds only to a small modification of the cost functional. 
This highlights another advantage of the CBO approach.

\subsection{Optimal control problems}
In this section we consider an optimal control problem constrained by a nonlinear nonlocal ODE, which models an interaction of sheep and shepherd dogs. Due to the nonconvexity induced by the nonlinear interactions, we expect to have local and global minima. The problem we discuss in the following was already considered in \cite{Burger2019} and a gradient descent (GD) algorithm was derived and used for numerical results. Of course, using a gradient approach we can only hope to convergence to local minima, while we aim for global solutions here.

\subsubsection{Mathematical formulation of the shepherd problem}
We now consider an optimal control problem that aims to gather the sheep close to a predefined destination point $x_\text{des} \in \mathbb R^2$ and therefore define the cost functional based on the variance and center of mass of the sheep flock. In more detail, let $N_S\in\mathbb N$ be the number of sheep in the flock and $x_i, v_i \colon [0,T] \rightarrow \R^2$ denote the position and velocity of the $i$-th sheep, respectively. We assume to be able to control the velocity of the dog $u \colon [0,T] \rightarrow \R^2$. The interaction of the sheep flock and the shepherd dog are assumed to be potential based, hence by Newton's Third Law we obtain the dynamics
\begin{subequations}\label{eq:sheepDog}
\begin{align}
&\frac{d}{dt}x_i = v_i, \qquad x_i(0)=x_{0,i}, \qquad i=1,\dots,N_S,\\ 
&\frac{d}{dt}v_i =-\alpha v_i -\frac{1}{N_S} \sum_{j=1}^{N_S} K_{s}\left(x_j - x_i\right) -  K_{d}\left(d-x_i\right), \qquad v_i(0)=v_{0,i}, \\ 
&\frac{d}{dt} d \,= u, \qquad d(0)=d_0,
\end{align}
\end{subequations}
where $K_s$ and $K_d$ are the gradients of the interaction potentials modeling binary sheep-sheep and sheep-dog interactions, respectively. The parameter $\alpha>0$ causes damping.

For the numerical implementation of the interactions we use the Morse potential \cite{DOrsogna2006}
\[
P\left(r\right) = C_{r_{s/d}}\exp\left(-\frac{|r|}{l_r}\right) - C_a \exp\left(-\frac{|r|}{l_a}\right), \quad K_{s/d} \left(r\right) = \left(-\frac{C_r}{l_r}\exp\left(-\frac{|r|}{l_r}\right) + \frac{C_a}{l_a} \exp\left(-\frac{|r|}{l_a}\right) \right) \frac{r}{|r|},
\]
where $ C_r, l_r $ are parameters for strength and range of the repulsion force, respectively, and $C_a, l_a$ are parameters for strength and range of the attraction force, respectively. We adjust the forces of the sheep-sheep and sheep-dog interactions by choosing different range and strength parameters. 

In order to model the task of gathering the sheep flock around the destination point $x_\text{des}$ we use the cost functional given by
\[
J = \int_0^T  \sigma_1 (V(t) - V_0)^2 + \sigma_2 \| E(t) - x_\text{des} \|^2 + \sigma_3\| u(t) \|^2  \, dt,
\]
where
$E(t) \in \mathbb{R}^D$ is the center of mass of the sheep flock at time $t$ defined as
\[
E(t) = \frac{1}{N_S} \sum_{i=1}^{N_S} x_i(t),
\]
and $V(t)$ is the variance of the sheep flock positions at time $t$ defined as
\[
V(t) = \frac{1}{N_S} \sum_{i=1}^{N_S} \| x_i(t) - E(t) \|^2.
\]
Hence, the first integrand in the cost functions measures the deviation of flock's variance from the prescribed target variance 
$V_0$. The second integrand measures the distance of the center of mass of the flock to the desired destination. Both together aim to gather the crowd around $x_\text{des}$. The third integrand has the physical meaning of minimizing the kinetic energy used for the task. Mathematically, it causes the cost functional to be coercive w.r.t.~$u \in L^2((0,T),\mathbb R^2)$. The weights $\lambda_1, \lambda_2,\lambda_3 \ge 0$ allow us to adjust the influence of the different terms. To have a regularity that is comparable to $L^2((0,T))$, which is used for the gradient descent approach, we choose $\nu=0.5$ for the regularity parameter of GP.

Assuming enough regularity for $K_{s/d}$, and given $u \in L^2((0,T))$ \eqref{eq:sheepDog} admits a unique solution $y(u) = (x_u,v_u,d_u)$. This allows us to define the reduced cost functional
$$f(u) = \mathcal J(y(u),u),$$
which is used in the CBO algorithm. We want to emphasize that for each evaluation of $f$ we have to solve the ODE system \eqref{eq:sheepDog}. We compare our results to the ones obtained by a gradient descent method  (GD). The latter only requires one solve of \eqref{eq:sheepDog} and its adjoint system in each gradient step. However, it only guarantees to find stationary points. 

\begin{figure}[t]
    \centering
    
    \begin{subfigure}[t]{0.48\textwidth}
        \centering
        \includegraphics[width=\linewidth]{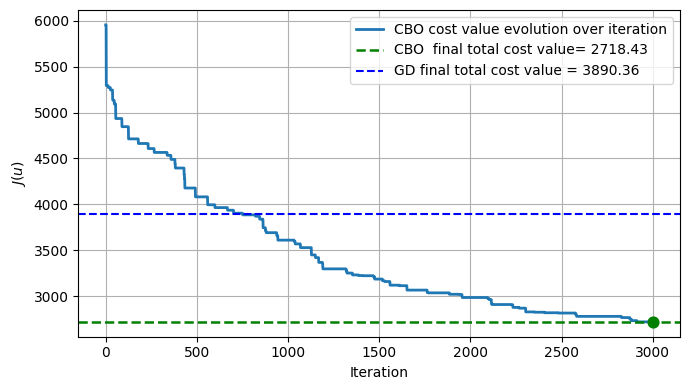}
        \caption{CBO cost evolution for a certain run in case of $\nu = 0.5$ with dashed reference lines indicating the final GD and CBO costs.}
        \label{fig:cbo_plot}
    \end{subfigure}
    \hfill
    \begin{subfigure}[t]{0.48\textwidth}
        \centering
        \includegraphics[width=\linewidth]{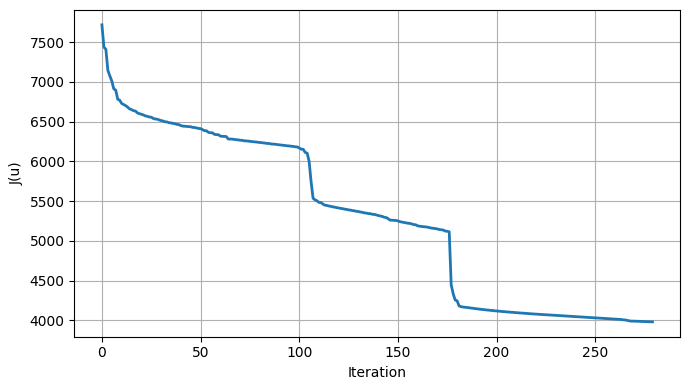}
        \caption{Cost evolution for the gradient descent (GD) method.}
        \label{fig:gd_plot}
    \end{subfigure}
    
    \caption{Comparison of the cost evolution over the iterations for CBO and GD.}
    \label{fig:cost_comparison}
\end{figure}

\begin{figure}[H]
    \centering
    \includegraphics[width=\linewidth]{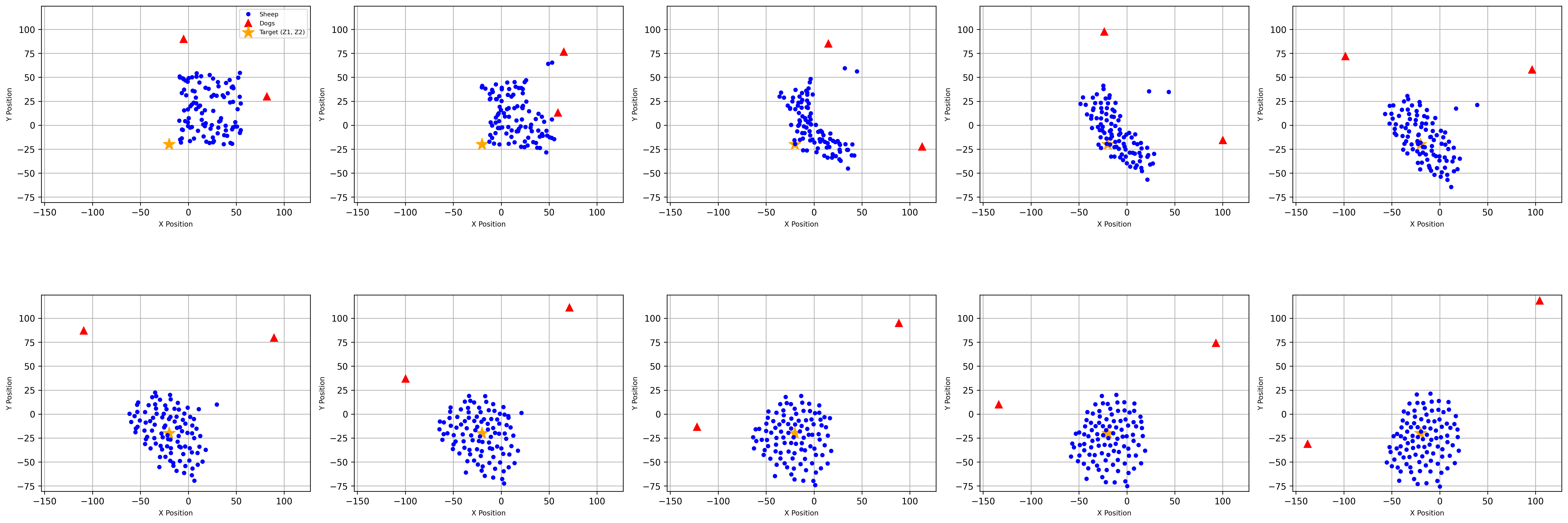}
    \caption{Snapshots of the sheep--dog dynamics under the optimized control for a certain run for value of $\nu = 0.5$. Blue circles represent the sheep flock, red triangles denote the dogs, and the orange star marks the target location.}
    \label{fig:sheep_dogs_evolution}
\end{figure}

In our specific case, the main optimization loop was executed for $3000$ iterations with $100$ CBO agents, 
resulting in a total of $3000 \times 100 = 300 000$ functional evaluations. 
This is considerably more expensive compared to the gradient descent method, which required only $2887$ evaluations for the step size search with Armijo rule and the update. Despite the higher computational cost of CBO compared to GD, CBO achieves a significant reduction of the objective functional. Indeed, we ran the CBO algorithm with $\nu=0.5$ and $\nu=2.5$ five times  each. For $\nu = 0.5$ the mean final total cost value is $2821.19$ and for $\nu = 2.5$ the mean final total cost value is $2840.94$ while the GD was stuck with a value of $3890.36$. According to the theory on GP, with $\nu=0.5$ we seek within functions of lower regularity than for $\nu=2.5$. Therefore the slight improvement for $\nu=0.5$ is reasonable.  The best run for $\nu = 0.5$ lead to nearly $30\%$ reduction in the cost compared to GD while the mean of the five runs of $\nu = 0.5$ lead to about $27.48\%$ reduction. The mean of the five runs of $\nu = 2.5$ lead to a reduction by $26.97\%$ compared to GD.
The cost evolution over the iterations for both methods is presented in Figure~\ref{fig:cost_comparison}. In the left panel the CBO cost evolution for a certain run in case of $\nu = 0.5$ is shown together with two dashed reference lines indicating the final total cost achieved by CBO and the final total cost obtained by GD. 
The right panel shows the cost evolution corresponding to the GD method.
The evolution of the sheep flock (blue circles) and the dogs (red triangles) under the optimized control for a value of $\nu = 0.5$ is shown in Figure~\ref{fig:sheep_dogs_evolution}. The simulation is performed with 100 sheep and 2 dogs, while the orange star indicates the target location.

\section{Conclusion and outlook}
In this article we have proposed a zero-order optimization algorithm for general problems in Sobolev spaces. The main idea was to exploit the toolbox of Gaussian Processes in order to sample functions with prescribed regularity that additionally satisfy constraints at given points. The algorithm was tested in different settings involving boundary value problems for ODE and PDE. Moreover, we have shown that the method is able to significantly improve the solution found by first-order methods for a highly nonlinear optimal control problem constrained by ODE.

The promising results provided in this paper motivate further research in various directions including the analysis of the proposed algorithm in the discrete setting, the convergence towards a nonlinear system of SDE in Sobolev Space as $\tau \to 0$. Moreover, it is open to check if the mean-field theory and in particular the global convergence proof that was established for CBO dynamics in finite state spaces can be generalized to the setting in Sobolev spaces as $N\to\infty$.

\section*{Acknowledgement} This research received funding by the German Federal Ministry of Education and Research under the grant number 05M22PXA. Responsibility for the content of this publication lies with the authors. 

\bibliographystyle{abbrv}
\bibliography{biblio}

\end{document}